\documentclass[reqno,centertags, 12pt]{amsart}
\usepackage{amsmath,amsthm,amscd,amssymb}
\usepackage{latexsym}
\newcommand{\bbR}{{\mathbb{R}}}
\newcommand{\bbD}{{\mathbb{D}}}
\newcommand{\bbZ}{{\mathbb{Z}}}
\newcommand{\bbC}{{\mathbb{C}}}

\newcommand{\calE}{{\mathcal E}}
\newcommand{\calH}{{\mathcal H}}
\newcommand{\calK}{{\mathcal K}}
\newcommand{\calL}{{\mathcal L}} 
\newcommand{\calM}{{\mathcal M}}

\newcommand{\dott}{\,\cdot\,}

\newcommand{\lb}{\label}
\newcommand{\f}{\frac}

\newcommand{\ti}{\tilde  }

\newcommand{\loc}{\text{\rm{loc}}}

\newcommand{\s}{\text{\rm{s}}}

\newcommand{\bi}{\bibitem}

\newcommand{\beq}{\begin{equation}}
\newcommand{\eeq}{\end{equation}}
\newcommand{\ba}{\begin{align}}
\newcommand{\ea}{\end{align}}

\allowdisplaybreaks
\numberwithin{equation}{section}

\newtheorem{theorem}{Theorem}[section]
\newtheorem*{t1}{Theorem 1}
\newtheorem*{t2}{Theorem 2}

\newtheorem{corollary}[theorem]{Corollary}
\theoremstyle{definition}

\theoremstyle{remark}

\newcommand{\abs}[1]{\lvert#1\rvert}

\begin{document}
\title{On a Theorem of Kac and Gilbert}
\author[B.~Simon]{Barry Simon}

\thanks{$^*$ Mathematics 253-37, California Institute of Technology, Pasadena, CA 91125. 
E-mail: bsimon@caltech.edu. Supported in part by NSF grant DMS-0140592} 

\date{April 16, 2004}

\begin{abstract} We prove a general operator theoretic result that asserts that many 
multiplicity two selfadjoint operators have simple singular spectrum.  
\end{abstract}

\maketitle

\section{Introduction} \lb{s1} 

In 1963, I.S.~Kac \cite{Kac} proved that whole-line Schr\"odinger operators, 
$-\f{d}{dx^2}+V(x)$,  for fairly general $V$'s have simple singular spectrum. 
It is well known (e.g., $V=0$) that the absolutely continuous spectrum can 
have multiplicity two and that under limit point hypotheses, eigenvalues 
are simple. But the simplicity of the singular continuous spectrum is 
surprisingly subtle. Some insight into the result was obtained by Gilbert 
\cite{Gil}, who found a proof using the subordinacy theory of Gilbert-Pearson 
\cite{GP}. The proof is elegant but depends on the substantial machinery of 
subordinacy. Our purpose here is to note an abstract result that relates 
these things to the celebrated result of Aronszajn-Donoghue \cite{AD}: 

\begin{t1} \lb{t1} Let $A$ be a bounded selfadjoint operator on $\calH$ and 
$\varphi\in\calH$ a cyclic vector for $A$. Suppose $\lambda\in\bbR
\backslash\{0\}$ and 
\begin{equation} \lb{1.1}
B=A+\lambda \langle\varphi,\dott\rangle \varphi 
\end{equation}
Then the singular spectral measures for $A$ and $B$ are disjoint. 
\end{t1} 

We state this and the next theorem in the bounded case for simplicity; we discuss the 
general case later. Here's the main result of this note:  

\begin{t2} \lb{t2} Let $\calH=\calK_1\oplus\calK_2$ and $P:\calH\to\calK_1$, the 
canonical projection. Let $A_j\in \calL(\calK_j)$ for $j=1,2,$ and $\varphi\in\calH$ so 
that $\varphi_1 =P\varphi$ and $\varphi_2 =(1-P)\varphi$ are cyclic for $A_1$ and 
$A_2$. Let $\lambda\in\bbR\backslash\{0\}$. Then 
\begin{equation} \lb{1.2}
C=A_1 \oplus A_2 +\lambda (\varphi, \dott)\varphi 
\end{equation} 
has simple singular spectrum. 
\end{t2} 

{\it Remark.} If $A_1$ is unitarily equivalent to $A_2$ and has a.c.~spectrum, then $C$ 
has multiplicity two a.c.~spectrum. So it is interesting that the singular spectrum is 
simple. 

\medskip
In Section~2, we prove Theorem~2. In Section~3, we apply it to whole-line Jacobi 
matrices. In Section~4, we discuss extensions of Theorem~2 to the case of unbounded 
selfadjoint operators and to unitary operators. In Section~5, we apply the results of 
Section~4 to Schr\"odinger operators and to extended CMV matrices.

\bigskip
\section{Proof of Theorem~2} \lb{s2} 

Let $d\mu_j$ be the spectral measure for $\varphi_j$ and $A_j$ and $B=A_1 \oplus A_2$. 
Thus $\calK_j \cong L^2 (\bbR,d\mu_j)$ in such a way that $A_j$ is multiplication by $x$ 
and $\varphi_j\cong 1$. 

Pick disjoint sets $X,Y,Z\subset\bbR$ whose union is $\bbR$ so $d\mu_1\restriction X$ is 
equivalent to $d\mu_2\restriction X$ and $\mu_1 (Y)=\mu_2(Z)=0$. For example, if 
$d\mu_1 =f(d\mu_1+ d\mu_2)$, then one can take $X=\{x\mid 0<f(x)<1\}$, $Y=\{x\mid 
f(x)=0\}$, $Z=\{x\mid f(x)=1\}$. 

Let $\calL_1$ be the cyclic subspace generated by $\varphi$ and $B$ and $\calL_2 = 
\calL_1^\perp$. Then $\psi =(\chi_X, -\chi_X)$ is a cyclic vector for $B\restriction 
\calL_2$ and its spectral measure is 
\begin{equation} \lb{2.1}
d\mu_\psi^B=\chi_X (x) (d\mu_1 + d\mu_2) 
\end{equation}
In particular, 
\begin{equation} \lb{2.2}
(d\mu_\psi^B)_\s \leq (d\mu_1 + d\mu_2)_\s 
\end{equation} 

By definition of $\calL_1$, $\varphi$ is cyclic for $B\restriction\calL_1$ and 
\[ 
C\restriction\calL_1 =B+\lambda (\varphi,\dott)\varphi 
\]
so, by Theorem~1, 
\begin{equation} \lb{2.3}
(d\mu_\varphi^C)_\s \perp (d\mu_\varphi^B)_\s = (d\mu_1 + d\mu_2)_\s 
\end{equation} 

Thus, the singular parts of $d\mu_\varphi^C$ and $d\mu_\psi^C=d\mu_\psi^B$ 
are disjoint, which implies that the singular spectrum of $C$ is simple. 
\qed 

\smallskip
The proof shows that the singular parts of $B$ and $C$ overlap in $\chi_X(d\mu_1 + 
d\mu_2)_\s$ and, in particular, 

\begin{corollary}\lb{C2.1} $B$ and $C$ have mutually singular parts if and only if 
$A_1$ and $A_2$ have mutually singular parts. 
\end{corollary} 

\bigskip
\section{Application to Jacobi Matrices} \lb{s3} 

A two-sided Jacobi matrix is defined by two two-sided sequences, $\{b_n\}_{n=-\infty}^\infty$ 
and $\{a_n\}_{n=-\infty}^\infty$ with $b_n\in\bbR$ and $a_n\in (0,\infty)$ and $\sup_n 
(\abs{a_n}+\abs{b_n})< \infty$. It defines a bounded operator $J$ on $\ell^2 (\bbZ)$ by 
\begin{equation} \lb{3.1}
(Ju)_n = a_{n-1} u_{n-1} + b_n u_n + a_n u_{n+1} 
\end{equation} 
so 
\begin{equation} \lb{3.2}
J=\begin{pmatrix} 
\cdots & \cdots & \cdots & \cdots & \cdots & \cdots\\
\ldots & b_{-1} & a_{-1} & 0 & 0 & \ldots \\  
\ldots & a_{-1} & b_0 & a_0 & 0 & \ldots \\
\ldots & 0 & a_0 & b_1 & a_1 & \ldots \\
\cdots & \cdots & \cdots & \cdots & \cdots & \cdots
\end{pmatrix}
\end{equation} 

\begin{theorem}\lb{T3.1} The singular spectrum of $J$ is simple. 
\end{theorem} 

\begin{proof} Let $\calK_1 = \ell^2 ((-\infty, -1])$, $\calK_2 =\ell^2 ([0,\infty))$, 
and $\varphi$ the vector with components 
\begin{equation} \lb{3.3}
\varphi_j = \begin{cases} 1 & j= -1,0 \\
0 & j\neq -1,0  \end{cases}
\end{equation} 
so $P\varphi = \delta_{-1}$; $(1-P)\varphi=\delta_0$. 

Then 
\[ 
J-a_{-1} (\varphi,\varphi)=A_1 \oplus A_2 
\]
where $A_2$ is the one-sided Jacobi matrix with 
\[
A_2 =\begin{pmatrix} 
b_0 -a_{-1} & a_0 & 0 & \dots \\
a_0 & b_1 & a_2 & \dots \\
\dots & \dots & \dots & \dots 
\end{pmatrix} 
\]
and $A_1$ in $\delta_{-1},\delta_{-2}, \dots$ basis is 
\[
A_1 =\begin{pmatrix} 
b_{-1} -a_{-1} & a_{-2} & 0 & \dots \\
a_{-2} & b_{-2} & a_{-3} & \dots \\
\dots & \dots & \dots & \dots 
\end{pmatrix} 
\]
Thus $P\varphi$ is cyclic for $A_1$ and $(1-P)\varphi$ is cyclic for $A_2$. 
Theorem~2 applies and implies the desired result. 
\end{proof}  

\bigskip
\section{Unitary and Unbounded Selfadjoint Operators} \lb{s4} 

Let $U_1$, $U_2$, and $W$ be unitary operators on $\calK_1$, $\calK_2$, and $\calH  
=\calK_1\oplus \calK_2$ so $W-U_1\oplus U_2$ is rank one, and so if $\varphi\in\ker 
(W-U_1\oplus U_2)^\perp$, then $P\varphi$ is cyclic for $U_1$ and $(1-P) 
\varphi$ is cyclic for $U_2$, where $P$ is the canonical projection of $\calH$ 
to $\calK_1$. Suppose $W-U_1\oplus U_2\neq 0$. Then  

\begin{theorem}\lb{T4.1} The singular spectrum of $W$\! is simple. 
\end{theorem} 

\begin{proof} We begin by proving the unitary analog of the Aronszajn-Donoghue 
theorem. If $W-V$ is rank one and nonzero, then for $\varphi$ a unit vector 
in $\ker (W-V)^\perp$, we have $W\varphi=\lambda V\varphi$ for some $\lambda 
\in\partial\bbD = \{z\in\bbC\mid \abs{z}=1\}$. By a direct calculation 
(see, e.g., \cite{OPUC}), 
\begin{align} 
\biggl(\varphi, \f{W+z}{W-z} \,\varphi\biggr) &= \f{1+zg(z)}{1-zg(z)} \lb{4.1} \\
\biggl( \varphi, \f{V+z}{V-z}\, \varphi\biggr) &= \f{1+zf(z)}{1-zf(z)} \lb{4.2}  
\end{align}
and 
\begin{equation} \lb{4.3}
g(z)=\lambda^{-1} f(z) 
\end{equation} 

If $\varphi$ cyclic for $W$\!, the singular spectrum in $W$ is supported on those 
$z\in\partial\bbD$ with 
\[
\lim_{r\uparrow 1}\, rz\, g(rz)=1
\]
and similarly, the singular spectrum of $V$ on the set of $z\in\partial\bbD$ with 
\[
\lim_{r\uparrow 1}\, rz\, f(rz)=1 
\]
By \eqref{4.3}, these sets are disjoint. 

This proves the Aronszajn-Donoghue theorem in the unitary case, and that implies 
this theorem by mimicking the proof of Theorem~2. 
\end{proof} 

Next, let $A_1$, $A_2$, and $C$ be potentially unbounded selfadjoint operators on 
$\calK_1$, $\calK_2$, and $\calH = \calK_1\oplus \calK_2$. Suppose $D\equiv 
(A_1\oplus A_2 -i)^{-1} -(C-i)^{-1}$ is rank one with $\varphi\in (\ker D)^\perp$ 
so $P\varphi$ is cyclic for $A_1$ and $(1-P)\varphi$ for $A_2$. Then with 
\[
U_j =(A_j+i)(A_j-i)^{-1} \qquad W=(C+i)(C-i)^{-1} 
\]
Theorem~\ref{T4.1} applies, so 

\begin{theorem}\lb{T4.2} $C$ has simple singular spectrum. 
\end{theorem}

\bigskip
\section{Extended CMV Matrices and Schr\"odinger Operators} \lb{s5} 

Extended CMV matrices enter in the theory of the orthogonal polynomials on the unit 
circle \cite{OPUC}. They are defined by a family of Verblunsky coefficients 
$\{\alpha_j\}_{j=-\infty}^\infty$ with $\alpha_j\in\bbD=\{z\in\bbC\mid \abs{z}=1\}$ 
as follows. Let $\Theta(\alpha)$ be the $2\times 2$ matrix 
\begin{equation} \lb{5.1}
\Theta (\alpha) = \left(\begin{array}{rr} 
\bar\alpha & \rho \\ 
\rho &-\alpha \end{array}\right) 
\end{equation}
where $\rho = (1-\abs{\alpha}^2)^{1/2}$. 

Think of $\ell^2(\bbZ)$, first as a direct sum $\oplus_{n=-\infty}^n \bbC^2$ with the 
$n$-th factor spanned by $(\delta_{2n},\delta_{2n+1})$ and let $\calM=\oplus\Theta 
(\alpha_{2n})$, then as a direct sum with the $n$-th factor spanned by $(\delta_{2n+1}, 
\delta_{2n+2})$ and $\calL=\oplus \Theta (\alpha_{2n+1})$. Then $\calE=\calL\calM$ is 
the extended CMV matrix. 

We claim 

\begin{theorem}\lb{T5.1} $\calE$ always has simple singular spectrum. 
\end{theorem} 

\begin{proof} Let 
\[
x=\f{1+\bar\alpha_{-1}}{1+\alpha_{-1}} 
\]
Then 
\[
\det \biggl( \Theta (\alpha_{-1}) - \begin{pmatrix} x & 0 \\ 0 & 1 \end{pmatrix}\biggr) 
=0 
\]
by a simple calculation. It is thus rank one, so if $\ti\calE$ is defined by replacing 
$\Theta (\alpha_{-1})$ by $\left(\begin{smallmatrix} x & 0 \\ 0 & 1 \end{smallmatrix} 
\right)$, $\ti\calE-\calE$ is rank one. On the other hand, $\ti\calE$ is a direct sum 
of two half-line CMV matrices and it is easy to see $P\varphi$ and $(1-P)\varphi$ are 
cyclic. Thus Theorem~\ref{T4.1} applies. 
\end{proof} 

Finally, we turn to the Schr\"odinger operator case that motivated us in the first place. 
Suppose $H=-\f{d^2}{dx^2} + V$ where $V\in L_\loc^1 (-\infty, \infty)$ is limit point at 
both $+\infty$ and $-\infty$. Let $H_1$ (resp.~$H_2$) be $H$ on $L^2 (0,\infty)$ 
(resp.~$L^2(-\infty,0)$) with $u(0)=0$ boundary conditions. Then 
\[
(H-i)- (H_1 \oplus H_2 -i)^{-1} 
\]
is rank one by the explicit Green's function formulae \cite{GS}, and its kernel is 
spanned by a function $\varphi$ with 
\begin{equation} \lb{5.2}
\biggl( -\f{d^2}{dx^2}+V\biggr) \varphi =i\varphi \qquad x\neq 0 
\end{equation}
with $\varphi$ $L^2$ at $+\infty$ and $-\infty$ and $\varphi (0_+)=\varphi(0_-)$. 
\eqref{5.2} implies $\varphi\restriction [0,\infty)$ (resp.~$\varphi\restriction (-\infty,0]$) 
is cyclic for $H_1$ (resp.~$H_2$). Theorem~\ref{T4.2} applies and yields the Kac-Gilbert 
theorem:  

\begin{theorem}\lb{T5.2} $-\f{d^2}{dx^2}+V$ has simple singular spectrum. 
\end{theorem}

\bigskip


\end{document}